\documentclass[12pt]{amsart}
\usepackage{amscd, amssymb}
\input{xy}
\xyoption{all}

\newtheorem{Thm}{Theorem}[subsection]

\newtheorem{Prop}[Thm]{Proposition}
\newtheorem{Def}[Thm]{Definition}
\newtheorem{Def/Thm}[Thm]{Definition/Theorem}

\theoremstyle{remark}

\newtheorem{EG}[Thm]{Example}

\numberwithin{equation}{subsection}
\newcommand{\ti }{\times}
\newcommand{\ot }{\otimes}
\newcommand{\ra }{\rightarrow}

\newcommand{\Proj}{{\mathrm{Proj}}}

\newcommand{\Spec}{{\mathrm{Spec}}}

\newcommand{\Ker}{{\mathrm{Ker}}}

\newcommand{\Pic}{{\mathrm{Pic}}}

\newcommand{\cA}{{\mathcal{A}}}
\newcommand{\cO}{{\mathcal{O}}}

\newcommand{\NN}{{\mathbb N}}

\newcommand{\PP }{{\mathbb P}}

\newcommand{\QQ }{{\mathbb Q}}
\newcommand{\CC }{{\mathbb C}}
\newcommand{\ZZ }{{\mathbb Z}}

\newcommand{\ka }{{\alpha}}

\newcommand{\kl }{{\lambda}}

\newcommand{\rk}{\mathrm{rk}}
\newcommand{\ext}{\mathrm{ext}}

\begin{document}

\title{Stable quasimaps}

\begin{abstract}
The moduli spaces of stable quasimaps unify various moduli appearing in the study of Gromov-Witten Theory.
This note is a survey article on the moduli of stable quasimaps,
based on papers \cite{CK, CKM, K} as well as the author's talk at Kinosaki Algebraic Geometry Symposium 2010.
\end{abstract}

\date{January 20, 2011, Revised April 18, 2011}

\author{Bumsig Kim}
\address{School of Mathematics, Korea Institute for Advanced Study,
Heogiro 87, Dongdaemun-gu, Seoul, 130-722, Korea}
\email{bumsig@kias.re.kr}

\maketitle
\section{Introduction}

A morphism from a variety $X$ to a projective space $\PP ^n$ is described by a linear system on $X$, which can
also be regarded as a $\CC ^\ti$-bundle $P$ on $X$ with a section $u$ of $P\times _{\CC ^\ti} \mathbb{A} ^{n+1}$ without base points.
When $X$ is a curve, one may compactify the morphism space
$\mathrm{Mor}(X, \PP ^n)$ by creating new rational components whenever
base points try to appear. This method eventually provides Kontsevich's stable map compactification.
There is another compactification, Quot scheme of rank $1$ subsheaves of $\mathcal{O}_X^{\oplus n+1}$.
The latter's boundary elements allow base points instead of attaching new rational components to $X$.
It turns out that this same idea can be applied to any GIT quotient $W/\!\!/G$ when a complex reductive Lie group
$G$ linearly acts on an affine variety $W$ with 
no strictly semistable points. (The precise definition and the condition for GIT quotients considered in this paper will be handled in
\S \ref{Quotients}.)

The above point of view leads us to the followings.

\begin{enumerate}

\item The notion of a quasimap $(P, u)$, i.e., a pair of a principal $G$-bundle $P$ on a prestable curve
and a section $u$ of $P\times _G W$ with at worst finitely many base points. 
The stability of quasimaps will be introduced in \S \ref{QNotion}.

\item New compactifications of moduli of  maps from curves
to a GIT quotient $W/\!\!/G$.   These include intermediate moduli spaces with moduli of
stable maps and moduli of stable quasimaps as the extremal ones on the parameter space of stabilities
(see  \cite{MM1, MM2, Toda} for $\PP ^n$ and Grassmannians;
the investigation of the general case
will appear elsewhere). In certain cases, the new spaces are easier to deal with than
the stable map spaces. We will present the stable quasimap compactification in \S \ref{quasimaps}.

\item The virtual smoothness of the Artin stack of the quasimap pairs when $W$ is LCI and $W/\!\!/G$ is smooth. The virtual smoothness will
be briefly treated in \S \ref{POT}.

\item  A new class
of examples with symmetric obstruction theory
if  $W/\!\!/G$ is a Nakajima quiver variety (see \cite{D, K}). This will be the topic of \S \ref{Hol}.

\item The wall-crossing interpretation of Givental's approach to the Classical Mirror Conjecture (see \cite{CK2}).

\end{enumerate}

\section{Three Quotients}\label{Quotients}

Let $W$ be an affine variety over $\CC$ with
a linear action by a complex reductive Lie group $G$. Typical examples of $G$
are products of general linear groups $GL_n (\CC )$.
In this situation one sometimes wants to define a quotient space.
There are three approaches.

\subsection{Affine quotients}
Since $W$ is affine, it can be considered as $\Spec$ of the ring $\CC [W]$ of all regular functions on $W$.
Hence, it is natural to define the quotient by $\Spec$ of the ring
$\CC[W]^G$ of all $G$-invariant regular functions on $W$. It is known that the projection
$W\ra \Spec\, \CC [W]^G$ is a good quotient of $W$ by the action of $G$ (see \cite[Theorem 6.3.1]{Le}).

In many cases, this is not interesting. For instance, if the homothety action is included in the $G$-action,
the closure of every $G$-orbit contains the origin. Therefore every $G$-invariant
function must be constant on $W$. Thus, $\Spec\, \CC[W]^G = \Spec\, \CC$. 

\subsection{GIT quotients}
The main reference for this subsection is King's paper \cite[\S 2]{King}.
In the previous example, to obtain an interesting space as a quotient,
we need to remove the origin. To do so geometrically, we should also prevent
other orbits from approaching the origin or some point.
For this, we will regard $W$ as $\Proj$ of the graded ring $\CC [W\times \mathbb{A}^1]$
whose grading comes from degrees with respect to the extra $\mathbb{A}^1$.
Given a character $\chi$ of $G$ (i.e., a one-dimensional representation $\chi$ of $G$),
we define a $G$-action on $W\times \mathbb{A}^1_{\chi ^{-1}}$, where
$\mathbb{A}^1_{\ka}$ for character $\ka$ denotes the one-dimensional representation space of $G$ associated to $\ka$.
Now it is natural to take $\Proj$ of the graded ring  $\CC [W\times \mathbb{A}^1]^G$ of all
$G$-invariant functions. This quotient is denoted by $W/\!\!/_{\chi} G$.

To describe the quotient space geometrically, let us recall the following notions (see \cite[Lemma 2.4]{King}).
 We call a point $p$ in $W$:
\begin{itemize}
\item {\em $\chi$-stable}  if,
for every one-parameter subgroup $\CC ^\ti$ of $G$, $ \CC ^\ti\cdot (p, 1)$  is  closed
in   $W\times \mathbb{A}^1$;

\item {\em $\chi$-unstable} if there is a one-parameter subgroup
$\CC ^\ti$ of $G$ such that the closure of $ \CC ^\ti\cdot (p, 1)$  meets $W\times \{0\}$;

\item {\em $\chi$-semistable} if it is not $\chi$-unstable.

\end{itemize}

\medskip 
Let $W^{ss}$ be the $\chi$-semistable locus of $W$ by the action of $G$.
If $f(x)\in \CC [W]$ and $f(x)z^l\in \CC[W\times \mathbb{A}^1]^G$ with $l > 0$,
then via the projection $W^{ss}\ra W/\!\!/_{\chi} G$ 
the inverse image of the open locus  $f(x)z^l \ne 0$ in $W/\!\!/_{\chi} G$ is $\Spec (\CC [W]_f)^G$.
This shows that $W/\!\!/_{\chi} G$ is a good quotient of $W^{ss}$.

For a suitable choice of $\chi$, it might be possible that there are no strictly semistable points. In such a case, the GIT quotient is
known to be also a geometric quotient so that $W/\!\!/_{\chi}G = W^s/G$ where $W^s$ is the stable locus, i.e., the locus of all $\chi$-stable points
(see \cite[\S 2]{King}).

\subsection{Stack quotients}
The stack quotient $[W/G]$, which is as a set (over $\Spec\, \CC$)
the set $W/G$ of $G$-orbits, keeps the data of isotropy subgroups (see \cite{LM}).
The stack quotient is defined as a category of groupoids over the category $\mathrm{Sch}$ of schemes over $\CC$,
whose objects over a scheme $S$
are pairs $(P, u)$ of principal $G$-bundles and $G$-equivariant morphism $u$ from $P$ to $W$. The morphism $u$
can be considered also as a section of $P\times _G W$.

\subsection{Relationships}
{\em Assume from now on that there are no strictly $\chi$-semistable points and $G$ acts on the stable locus
freely.} The latter condition is only technical (see \cite{C}).

With this assumption the three quotients are related by the diagram
\[\xymatrix{ [W/G ] \ar@{<-^{)}}[r]_{\text{open}} & W/\!\!/_{\chi}G \ar[rr]_{\text{projective}} & &\Spec\, \CC [W]^G }.\]
(Note that $W/\!\!/_{\chi =0} G$ coincides with the affine quotient $\Spec\, \CC [W]^G $, but we will not use this notation.)

\subsection{Examples}
There are many examples.

\begin{EG} Let $Y$ be a projective variety in $\PP ^n$ and let $C(Y)$ denote the affine cone of $Y$ in $\mathbb{A}^{n+1}$.
Then $Y$ is the GIT quotient $C(Y)/\!\!/_{\mathrm{det}}\CC ^\ti$, where the character is defined by the determinant map $\mathrm{det}$.
Typical examples for this case will be smooth complete intersections $Y$.
\end{EG}

\begin{EG} The smooth toric varieties can be considered as a GIT quotient $Y=\CC ^N /\!\!/_{\chi} G$ (see \cite{Au, Cox, Dol}). 
Since $G=\mathrm{Hom}(A_{n-1}(Y),\CC ^\times )$, we may regard an element in $A_{n-1}(X)$ as a character.
Then any very ample line bundle can be the character $\chi$.    

The partial flag varieties of nested subspaces of dimension $k_1, ..., k_r$ in $\CC ^{n=k_{r+1}}$
can be considered as a GIT quotient $$\prod _{i=1}^r \mathrm{Hom}(\CC ^{k_i},\CC ^{k_{i+1}})/\!\!/_{\chi}\prod _{i=1}^r GL_{k_i}(\CC )$$ with 
$\chi$ being the product of determinant characters of each components (for instance, see \cite{BCK}). Hence,
complete intersections $Y$ in smooth toric varieties  or partial flag varieites are also such GIT examples.
\end{EG}

\begin{EG} Quiver varieties are such GIT examples (for instance, see \cite{Gin}). 
In particular, quiver varieties of Nakajima type, being noncompact, will be recalled in section \ref{Hol}.
\end{EG}

\bigskip

\section{Quasimaps}\label{quasimaps} 

A morphism from $X$ to the GIT quotient $W^s/G$
amounts to a pair $(P, u)$ where $P$ is a principal $G$-bundle on $X$ and $u$ is a $G$-equivariant
map from $P$ to $W$ whose image is contained in $W^s$. By Luna's slice theorem,
$W^s$ is a principal $G$-bundle on $W^s/G$ in \'etale topology.
The direction $\Rightarrow$ therefore follows.
The other direction holds since $P\ra X$ is a categorical quotient.
Here we exchange left and right actions via the
inverse map $G\ra G, g\mapsto g^{-1}$.

\subsection{Stable quasimaps}\label{QNotion} If we allow that $u$ hits the unstable locus $W^{un}$,
the pair $(P, u)$ is not any more a well-defined map to the GIT quotient, but it
is a map $[u]$ to the stack quotient by the very definition of the stack quotient.

\begin{Def} The pair $(P, u)$ is called a {\em quasimap} of genus $g$ if:
\begin{itemize}
\item $X$ is a projective, smooth or nodal curve of genus $g$, with $n$ ordered smooth distinct markings;
\item the base locus $u^{-1}(P\times _G W^{un})$ consists of finite points.
\end{itemize}

The quasimap is called a {\em stable quasimap} if:
\begin{itemize}
\item $\omega _X\ot (P\times _G \mathbb{A}^1 _\chi )^\epsilon$ is ample for all positive rational number $\epsilon$;
\item the base points (if any) are smooth and non-marked points.
\end{itemize}
\end{Def}

\subsection{Degrees} Note that there is a notion of degree for a map $f$ from curve $X$ to $W/\!\!/G$ by the homomorphism
\[ \begin{array}{rcc} \deg (f) : \Pic (W/\!\!/_{\chi}G) &\ra& \ZZ \\
                                          M &\mapsto & \deg (f^*M) \end{array} .\]

Similarly we define the degree of $(P, u)$ as a homomorphism from the character group of $G$
to $\ZZ$ by sending $\ka\mapsto \deg  (P\times _G \mathbb{A}^1 _{\ka})$.
Let $\beta$ be
a group homomorphism from the character group of $G$ to $\ZZ$.

\begin{Thm} {\em (\cite{CK, CKM})}
The moduli stack $Q_{g, n} (W/\!\!/_{\chi}G, \beta )$ of stable quasimaps
of $n$-pointed, genus $g$, degree $\beta$ to $W/\!\!/_{\chi} G$ is a finite-type DM
stack  proper over the affine quotient $\Spec\, \CC [W] ^G$.
Furthermore if the affine scheme $W^s$ is smooth and $W$ is a locally
complete intersection scheme, then the moduli stack comes with a natural perfect obstruction theory.
\end{Thm}

The map from the moduli stack to the affine quotient can be naturally given since $P\ra X$ is categorical at the diagram
\[\xymatrix{ P \ar[r] \ar[d]\ar[dr] & W \ar[d] \\
X \ar@{.>}[r] & \Spec \CC [W] ^G } \] and $X$ is a projective scheme over $\CC$.

\subsection{Completeness} Let $\Delta$ be a unit disk in complex line $\mathbb{A}^1$ and let $\Delta ^\circ = \Delta \setminus \{0\}$.
Fix a smooth curve $X$ and consider $\Delta ^\circ$-family $(P, u)$ of stable quasimaps on 
$X\times \Delta ^\circ$. Suppose that $G=GL_r(\CC )$ so that
we may have the corresponding locally free sheaf $\mathcal{E}$ associated to the $P\times _{GL_r (\CC ) } \CC ^r$. In this situration, we 
indicate how to find an extension of $(P, u)$.
First we take a coherent sheaf extension $\overline{\mathcal{E}}$ of $\mathcal{E}$ using  the properness of Quot schemes (after showing the 
boundedness). 
The double dual $\overline{\mathcal{E}}^{\vee\vee}$
is reflexive on the smooth complex surface $C\times \Delta$ being regular in codimension 3 (see \cite[Lemma 1.1.10]{OSS}),
hence it is a locally free sheaf on the surface. 
Let $\overline{P}$ be the associated principal 
$G$-bundle on $X\times \Delta$.
Now we may consider $u$ as a rational map  \xymatrix{ C\times \Delta  \ar@{-->}[r] & \overline{P}\times _G W }.
The indeterminant locus of $u$ is a locus of finite points. However by Hartogs' theorem we conclude that the locus
is empty.

\subsection{Perfect obstruction theory}\label{POT}  It is not difficult to see
that the deformations of curves $X$ and principal bundles $P$ on $X$
are unobstructed. Let us fix $X$, $P$ and deform only section $u$ of $P\times _G W$. Note that the deformation space $\mathrm{def} (u)$ is
$\mathbb{H}^0 (X, u^*{\bf T}_\rho )$ where ${\bf T}_\rho $ is the relative tangent complex of the natural projection $\rho :P\times _G W \ra X$.
We can then show that $\mathrm{ob} (u) :=\mathbb{H}^1(X, u^*{\bf T}_\rho )$ is an obstruction space.
The relative virtual dimension is $\dim \mathrm{def} (u) - \dim \mathrm{ob} (u)$, which is locally constant.
Roughly speaking, this implies that
the moduli space is virtually smooth in the following sense. 
There is a vector bundle stack $[F/E]$ of a two term perfect complex $E\ra F$ whose kernel is the deformation sheaf and whose cokernel is
the obstruction sheaf and there is a closed embedding of
an intrinsic normal cone $C_Q$ of the moduli stack $Q=Q_{g, n }(W/\!\!/G, \beta )$ so that the refined intersection
$0^!_{F/E} (C_Q)$ of $C_Q$ with the zero section is defined to be the so-called virtual fundamental class 
$[Q_{g,n}(W/\!\!/G,\beta )]^{\mathrm{vir}}\in
A _{\mathrm{exp. dim}}(Q_{g,n}(W/\!\!/G,\beta ), \QQ)$ where the expected dimension $\mathrm{exp.dim}$ is the relative virtual dimension plus the dimension
of the $\mathfrak{M}_{g,n}$-relative stack $\mathfrak{B}un_G$ of principal $G$-bundles on the universal curve of the stack $\mathfrak{M}_{g,n}$ of $n$-pointed
genus $g$ prestable curves (see \cite{BF, Kr, LT}).

\subsection{Historical remarks}

These are limited remarks.

\subsubsection{} For a projective smooth toric variety, the spaces of stable quasimaps with the fixed domain curve $\PP ^1$
also become projective smooth toric varieties. The spaces are used to prove the Mirror Theorem for Fano/CY complete intersections by Givental in \cite{G}.

\subsubsection{}
Let $W/\!\!/_{\chi} G = \mathrm{Hom}(\CC ^r, \CC ^n)/\!\!/_{\det} GL_r( \CC) = \mathrm{Grass}(r, n)$. In this case,
a stable quasimap amounts to a rank $n-r$ quotient of $\cO _X ^{\oplus n}$ on a prestable curve with
certain conditions: for example when there are no
markings, the conditions are no torsion at nodes and no rational tails.
The latter is called a stable quotient and introduced
by Marian, Oprea, and Pandharipande in \cite{MOP}.

For a fixed smooth curve $X$, the moduli spaces of stable quotients are nothing but
Quot schemes of rank $n-r$ quotients of  $\cO _X ^{\oplus n}$.
Quot schemes have been used and studied in Gromov-Witten theory (for instance, see \cite{Ber, OT}).

\subsubsection{}
For a smooth projective toric variety $\mathbb{C} ^N/\!\!/ (\CC ^\times )^r$, the theorem is proven in \cite{CK}.
The paper \cite{CK} shows the idea that all the above constructions can be unified and generalized to any GIT
quotient $W/\!\!/G$ considered in section \ref{Quotients}.

\subsection{Quasimap invariants}
Using the perfect obstruction theory on the moduli space of stable quasimaps,
we define the virtual fundamental class of the moduli space.
 Hence, we can define intersection numbers
by integrals of tautological cohomology classes against the virtual fundamental class.

We conjecture that {\em these invariants and Gromov-Witten invariants for $W/\!\!/G$ carry the same amount of information} (see \cite{CK, CKM, CK2}).
A precise formulation of the conjecture is unknown except for the following two cases.

\begin{enumerate}
\item When $W/\!\!/G$ has the property that for every curve $C$
in $W/\!\!/G$, $C\cdot K_{W/\!\!/G} \le -2$, we expect that
both invariants exactly coincide.

\item The genus zero quasimap invariants should lie on the Lagrangian cone generated by the genus zero
gravitational Gromov-Witten invariants (and vice-versa).
\end{enumerate}

\subsection{Some evidence}

\subsubsection{} For {\em fixed} $X=\PP ^1$ and any toric complete intersection $W/\!\!/G$, (2) is the Mirror Theorem in \cite{G}.
For Grassmannian case, (1) is a theorem  in \cite{MOP}.

\subsubsection{} In \cite{CK2} we prove (2) for any Fano/CY toric complete intersection in a toric variety and
(1) for any Fano toric variety.

\section{Stable Quasimaps to Holomorphic Symplectic Quotients}\label{Hol}

Let $V$ be a smooth affine variety with a holomorphic symplectic form $\omega$. The form
$\omega\in \Gamma (V, \mathcal{T}_V)$ is a nondegenerate closed (2,0)-form. In this setup one can define
a quotient which is also holomorphic symplectic.

\subsection{Holomorphic symplectic quotients} Suppose that the $G$-action $V$ is hamiltonian which means that:
the $G$-action preserves $\omega$
and there is a $G$-equivariant morphism $\mu :V \ra \mathfrak{g} ^*$ such that
\begin{equation}\label{duality} \tag{$\star$} \langle d\mu (\xi), g\rangle = \omega (\xi , d\alpha (g ) ),  \xi \in T_V. \end{equation}
Here $\ka :G \ra \mathrm{Aut} V$ is the map induced from the action and $\mathfrak{g}$ denotes the Lie algebra of
$G$.
The morphism $\mu$ is called a complex moment map.

We define the holomorphic symplectic quotient by $\mu ^{-1} (\kl ) /\!\!/ _{\chi } G$
where $\kl$ is a $G$-invariant regular value of
$\mu$. The quotient is denoted by $V/\!\!/\!\!/_{\kl, \chi } G$.

\subsection{Self-duality}
Let $\bf F$ denote the $G$-equivariant complex
$$[\mathfrak{g}\ot\mathcal{O}_V \stackrel{d\ka}{\ra}
\mathcal{T}_V \stackrel{d\mu}{\ra} \mathfrak{g}^*\ot\mathcal{O}_V ] _{|_{\mu ^{-1}(\kl )}}.$$
Note that by \eqref{duality} the self-duality ${\bf F}\cong {\bf F}^\vee$ holds.
Note also that ${\bf F}$ is a monad such that
$\Ker d\mu / \mathrm{Im} d\ka  $ in $(\mu ^{-1} (\kl ))^s$ is
 isomorphic to the pullback of the tangent sheaf of holomorphic
symplectic quotient. In other words, the quotient ${\bf F}/G$ 
is a generalized Euler sequence
for the tangent sheaf of $V\!\!/\!\!/\!\!/G$.

\subsection{Symmetric obstruction theory}

Let us fix a smooth projective curve $X$.
We define $\mathfrak{M}_\beta$ to be the stack of degree $\beta$ stable quasimaps
to the holomorphic symplectic quotient $\mu ^{-1}(\kl) /\!\!/G$ from $X$.
This moduli space has a symmetric obstruction theory if $X$ is an elliptic curve, i.e.,
the deformation space $\mathbb{H}^0(X, P\times _G {\bf F})$ is functorially isomorphic to the dual of
the obstruction space $\mathbb{H}^1(X, P\times _G {\bf F})$. This follows from the Serre duality and the self-duality ${\bf F}\cong {\bf F}^\vee$.

By twisting, this symmetry can be made hold for arbitrary genus, smooth projective curve $X$ 
when the quotient is a Nakajima quiver variety.

\subsection{Nakajima's quiver varieties}

 A quiver $Q$ is an oriented graph, i.e.,
 data $(Q_0, Q_1, h, t)$ where $Q_0$ is the set of vertices, $Q_1$ is the set of arrows,
 $h$ is the head map, and $t$ is the tail map ($h, t: Q_1 \ra Q_0$).
 Let $\overline{Q}$ be the quiver obtained from $Q$ by adding
 the opposite arrow $\bar{a}$ for each arrow $a$ in $Q_1$ (so $|\overline{Q}_1| = 2 |Q_1|$).
 Set $\bar{\bar{a}}=a$.

Fix a distinguished vertex $0\in Q_0$ and a dimension vector $v\in \NN ^{Q_0}$.
Let $V$ be the direct sum of $\mathrm{Hom} (\mathbb{C}  ^{v_{ta}}, \mathbb{C}  ^{v_{ha}})$ for all $a\in \overline{Q}_1$.
$V$ has a decomposition $V_+\oplus V_-$ based on the arrows in $Q_1$ and $\overline{Q}_1\setminus Q_1$
so that it is a symplectic vector space with the canonical symplectic form.
Also $V$ comes with a natural action of $G=\Pi _{i\in Q_0, i\ne 0} GL _{v_i} (\CC )$.
It is easy to see that this action is hamiltonian with a moment map
$$\left(\sum _{a\in \overline{Q}_1: ha =i} (-1)^{|a|} \phi _a \circ \phi _{\bar{a}}\right) _{i \in Q_0, i\ne 0}.$$
Choose $(\theta _i)\in \ZZ ^{Q_0\setminus \{0\}}$ and define $\theta :G\ra \CC ^\ti $ by $g\mapsto \prod \det g_i ^{\theta _i}$.
Now we can define $V/\!\!/\!\!/_{\lambda ,\theta} G$. Let $\lambda =0$ from now on.

Note that a quasimap data
is equivalent to
$(E_i, \phi _a )$ where $E_i$ is a vector bundle $P\times _G \mathbb{A} ^{v_i}$ and
$\phi _a : E_{ta}\ra E_{ha}$ is a homomorphism obtained from $u$, satisfying the moment map relation
and the condition that
\lq base points' are finite.
This therefore motivates the following.

\subsection{Twisted quiver sheaves}

Fix a smooth projective curve $X$,
a line bundle $M_a$ on $X$ for every $a\in \overline{Q}_1$, and
an isomorphism $M_a\ot M_{\bar{a}} \ra K_X^{-1}$ for every $a\in Q_1$.

\begin{Def} $(E_i, \phi _a)$ is called a {\em framed-twisted quiver sheaf on $X$ with the moment map relation} if:
\begin{enumerate}
\item $E_i$ is a coherent sheaf on $X$;

\item $\phi _a : M_a\ot E_{ta}  \ra E_{ha}$ is an $\mathcal{O}_X$-homomorphism;

\item $\sum _{i\ne 0} \sum _{t\bar{a}=i} (-1)^{|a|} \phi _a \circ (\mathrm{Id}_{M_{a}}\ot\phi _{\bar{a}})= 0$ in
$\bigoplus _{i\ne 0} \mathrm{End}(K_X^{-1}\ot E_i, E_i)$;

\item $E_0 = \mathcal{O}_X^{\oplus r}$ for some integer $r$.

\end{enumerate}
\end{Def}

Often, we will call it simply twisted quiver sheaf, even just quiver sheaf.

\subsection{Stability Conditions}

Let $\mathcal{A}$ be the abelian category of twisted quiver sheaves with respect to the fixed data above.
Let $Z : K(\cA)\ra \CC$ be a homomorphism defined by
$$ (E_i, \phi _a )\mapsto \sum _{i\ne 0}\rk E_i + \sqrt{-1} (\sum _{i\ne 0} \deg E_i + \tau \rk E_0), $$
where $\tau$ is a positive real number.
Note that $Z(K(\cA )\setminus \{0\})$ is contained in the union of the half plane where the real part is positive
with the positive $y$-axis.
This is a Bridgeland stability function with Harder-Narasimhan property.

A nonzero twisted quiver sheaf $(E_i, \phi _a)$ is called $\tau$-(semi)stable
if $\mathrm{Arg} Z(E') < (\le) \mathrm{Arg} Z(E)$
for every nonzero proper subobject $E'$ of $E$ in $\cA$.

\begin{Prop} {\em (\cite{K})}
Fix a dimension vector $v$  with $v_0=1$ and a degree vector $d$.
Then there is a positive number $\tau _0$ such that for every $\tau\ge \tau _0$ and for every
$(E_i, \phi _a)$ twisted quiver sheaf with $(v,d)$, TFAE

\begin{enumerate}
\item $\tau$-semistability.

\item $\tau$-stability.

\item Stability as a twisted quasimap to $V/\!\!/\!\!/G$ with $\theta = (1,...,1)$.

\end{enumerate}

\end{Prop}

Here the stability as a twisted quasimap $(P, u)$ is defined to be similar to the untwisted case (for the detail see \cite{K}): 
After fixing a local trivialization of $M_a$ and $\omega _X$,
we may say whether a point $X$ lands via $u$ on the $\theta$-unstable locus of $W$ or not. The stability is by definition 
the requirement that the generic point of $X$ lands on
the stable locus $W^s = W^{ss}$. This condition is independent of the choices of local trivializations. 

\medskip

Fix $(v,d)$ with $v_0=1$. For $\tau$ to be a wall (i.e., there is strictly $\tau$-semistable objects with type $(v,d)$), 
$\tau$ must satisfy a numerical condition that the $\tau$-slope of $(v,d)$ equals to
the $\tau$-slope of some $(v', d')\ne (v,d)$ with $v'_i\le v_i$. 
This condition ensures that the walls for a given $(v,d)$ are discrete, which together with the above 
proposition shows that 
there are only finitely many walls $\tau _i$, $i=1,...,N$.

\begin{Thm} {\em (\cite{K})}
The moduli stack  $\mathfrak{M}^\tau _{(v,d)}$ of $\tau$-stable quiver
sheaves of the rank-degree vector $(v,d)$ is a finite type algebraic space equipped with
a symmetric obstruction theory. The stack is proper over $\Spec \, \mu ^{-1}(0) ^G$ if
$\tau$ is large enough.
\end{Thm}

\subsection{Remarks}
When $\overline{Q}$ is the ADHM quiver, the above proposition and theorem were proven in \cite{D} which
is one of the main sources of inspiration for the general case.

So far, there are two classes of moduli examples equipped with symmetric obstruction theory:
moduli of stable objects in the abelian category of coherent sheaves of a CY 3-fold
and representations of a quiver with relations from a superpotential on $Q$
(see \cite{Th, PT, S, JS, KS}).

One can study wall-crossings of topological Euler characteristics of  $\mathfrak{M}^\tau _{(v,d)}$ weighted by Behrend's
constructible functions (see \cite{Behrend}) using Joyce-Song formula (see \cite{JS}) or Kontsevich-Soibelman formula (see \cite{KS}).
Once again for the ADHM case, it is done in
\cite{D, CDP1, CDP2}; and its generalization is carried out in \cite{KL}.

For the use of Joyce-Song theory, in \cite{KL} we prove:

\begin{itemize}

\item $\chi (E, F) : = \ext ^0 _{\mathcal{A}}(E, F) - \ext ^1 _{\mathcal{A}}(E, F) + \ext ^1 _{\mathcal{A}}(F,E) - \ext ^0 _{\mathcal{A}}(F, E)$
is numerical in when $E, F \in \mathcal{A}$ are locally free quiver sheaves with $\rk E_0 + \rk F_0\le 1$;

\item $\mathfrak{M}^\tau _{(v,d)}$ is analytic-locally a critical locus of a holomorphic function
on a smooth analytic domain.

\end{itemize}

It would be interesting to relate $\cA$ with a CY 3-category.

\subsection{Acknowledgments} The author thanks Ionu\c{t} Ciocan-Fontanine, Changzheng Li, and the Referee for invaluable
 comments on this survey article. This work is financially supported by KRF-2007-341-C00006.


\begin{thebibliography}{9999}

\bibitem{Au} M. Audin,  {\em The topology of torus actions on symplectic manifolds,}
 Progress in Mathematics, 93. Birkh\"auser Verlag, Basel, 1991.

\bibitem{Behrend} K. Behrend, {\em Donaldson-Thomas type invariants via microlocal geometry,}
Ann. of Math. (2) 170 (2009), no. 3, 1307-1338.

\bibitem{BF} K. Behrend, B. Fantechi,  {\em The intrinsic normal cone,}
Invent. Math. 128 (1997), no. 1, 45-88. 

\bibitem{Ber} A. Bertram, {\em Quantum Schubert calculus,} Adv. Math. 128 (1997), no. 2, 289-305.

\bibitem{BCK} A. Bertram, I. Ciocan-Fontanine, B. Kim {\em Gromov-Witten invariants for abelian and nonabelian quotients,}
 J. Algebraic Geom. 17 (2008), no. 2, 275-294.

\bibitem{C} D. Cheong, in preparation.

\bibitem{Cox} D. Cox,  {\em The homogeneous coordinate ring of a toric variety,}
J. Algebraic Geom. 4 (1995), no. 1, 17-50. 

\bibitem{CDP1} W-E. Chuang, D.E. Diaconescu, G. Pan, {\em Chamber structure and wall- crossing in the ADHM theory of curves II},
         {\tt arXiv:0908.1119}.
\bibitem{CDP2}  W-E. Chuang, D.E. Diaconescu, G. Pan, {\em Rank Two ADHM Invariants
and Wallcrossing,} {\tt arXiv:1002.0579}.

\bibitem{CK} I. Ciocan-Fontanine, B. Kim, {\em Stable toric quasimaps}, 
Advances in Mathematics, Volume 225, Issue 6, 20 December 2010, 3022-3051.

\bibitem{CK2} I. Ciocan-Fontanine, B. Kim, in preparation.

\bibitem{CKM} I. Ciocan-Fontanine, B. Kim, D. Maulik, {\em Stable quasimaps to GIT quotients,} in preparation.


\bibitem{D} D.-E. Diaconescu, {\em Chamber structure and wallcrossing in the ADHM theory of
curves I,} {\tt arXive:0904.4451.}

\bibitem{Dol} I. Dolgachev, {\em Lectures on invariant theory,}
 London Mathematical Society Lecture Note Series, 296. Cambridge University Press, Cambridge, 2003.

\bibitem{Gin} V. Ginzburg, {\em Lectures on Nakajima's quiver varieties,} 	{\tt arXiv:0905.0686}.

\bibitem{G} A. Givental, {\em A mirror theorem for toric complete intersections}, in \lq\lq Topological field theory, primitive
forms and related topics (Kyoto, 1996)", Progr. Math., 160, Birkh\"auser Boston, Boston, MA, 1998, 141--175.

\bibitem{JS} D. Joyce, Y. Song, {\em
A theory of generalized Donaldson-Thomas invariants,}
    to appear in Memoirs of the A.M.S, {\tt arXiv:0810.5645}.

\bibitem{K} B. Kim, {\em Stable quasimaps to holomorphic symplectic quotients},  {\tt arXiv:1005.4125}.

\bibitem{KL} B. Kim,  H. Lee, {\em Wall-crossings for twisted quiver bundles,} {\tt arXiv:1101.4156}.

\bibitem{King} A.D. King, {\em Moduli of representations of finite-dimensional algebras,}
 Quart. J. Math. Oxford Ser. (2) 45 (1994), no. 180, 515-530.

\bibitem{KS}  M. Kontsevich, Y. Soibelman,
{\em Stability structures, motivic Donaldson-Thomas invariants and cluster transformations,} {\tt arXiv:0811.2435}.

\bibitem{Kr} A. Kresch , {\em Cycle groups for Artin stacks,}
Invent. Math. 138 (1999), no. 3, 495-536.

\bibitem{LM} G. Laumon, L. Moret-Bailly, {\em Champs alg\'ebriques,}
A Series of Modern Surveys in Mathematics, 39. Springer-Verlag, Berlin, 2000.

\bibitem{Le} J. Lepotier, {\em Lectures on vector bundles,}
Translated by A. Maciocia. Cambridge Studies in Advanced Mathematics, 54. Cambridge University Press, Cambridge, 1997.

\bibitem{LT} J. Li, G. Tian, {\em Virtual moduli cycles and Gromov-Witten invariants of algebraic varieties,}
 J. Amer. Math. Soc. 11 (1998), no. 1, 119-174.

\bibitem{MOP} A. Marian, D. Oprea,  R. Pandharipande, {\em The moduli space of stable quotients,} {\tt arXiv:0904.2992}.


\bibitem{MM1} A. Musta\c{t}\v{a}, A. Musta\c{t}\v{a}, {\em Intermediate moduli spaces of stable maps,} Invent. Math. 167 (2007), no. 1, 47-90.

\bibitem{MM2} A. Musta\c{t}\v{a}, A. Musta\c{t}\v{a}, {\em
The Chow ring of $\overline{M}_{0,m}(\PP ^n, d)$,} J. Reine Angew. Math. 615 (2008), 93-119.


\bibitem{OSS} C. Okonek, M. Schneider,  H. Spindler, {\em Vector bundles on complex projective spaces,}
 Progress in Mathematics, 3. Birkh\"auser, Boston, Mass., 1980.

\bibitem{OT} C. Okonek, A. Teleman, {\em Comparing virtual fundamental classes: gauge theoretical Gromov-Witten invariants for
toric varieties,} Asian J. Math. 7 (2) (2003), 167-198.

\bibitem{PT} R. Pandharipande, R.P. Thomas, {\em Curve counting via stable pairs in the
derived category,} Invent. Math. 178 (2009), no. 2, 407-447.


\bibitem{S} B. Szendr\"oi, {\em Non-commutative Donaldson-Thomas invariants and the conifold,} Geom. Topol. 12 (2008), no. 2, 1171-1202.

\bibitem{Th} R.P. Thomas, {\em A holomorphic Casson invariant for Calabi-Yau 3-folds, and
bundles on K3 fibrations,} J. Differential Geom. 54 (2000), no. 2, 367-438.

\bibitem{Toda} Y. Toda, {\em  Moduli spaces of stable quotients and the wall-crossing phenomena,}  {\tt arXiv:1005.3743}.


\end{thebibliography}
\end{document}